\documentclass[10pt]{article}
\usepackage{mathtools}
\usepackage{amsmath}
\usepackage{amssymb}
\usepackage{braket}
\usepackage{amsmath,color}
\newtheorem{theorem}{Theorem}[section]

\newtheorem{proposition}[theorem]{Proposition}
\newtheorem{definition}{Definition}[section]

\newtheorem{hypothesis}[theorem]{Hypothesis}

\newtheorem{remark}[theorem]{Remark}

\newcommand{\si}{\sigma}

\newcommand{\de}{\delta}
\newcommand{\De}{\Delta}
\newcommand{\om}{\omega}
\newcommand{\Om}{\Omega}

\newcommand{\f}{\mathcal{F}}
\newcommand{\g}{\mathcal{G}}

\newcommand{\p}{\mathcal{P}}

\newcommand{\bee}{\begin{equation}}
\newcommand{\eee}{\end{equation}}
\newcommand{\bea}{\begin{eqnarray}}
\newcommand{\eea}{\end{eqnarray}}
\newcommand{\bean}{\begin{eqnarray*}}
\newcommand{\eean}{\end{eqnarray*}}
\makeatletter
\@addtoreset{equation}{section}

\makeatother
\def\qed{{\hfill\hbox{\enspace${ \square}$}} \smallskip}
\def\sqr#1#2{{\vcenter{\vbox{\hrule height .#2pt \hbox{\vrule
 width .#2pt height#1pt \kern#1pt \vrule
width .#2pt} \hrule height .#2pt}}}}
\def\square{\mathchoice\sqr54\sqr54\sqr{4.1}3\sqr{3.5}3}
\def\ds{\begin{displaystyle}}
\def\eds{\end{displaystyle}}
\def\dis{\displaystyle }
\def\<{\langle }
\def\>{\rangle }

\def\R{\mathbb R}

\def\E{\mathbb E}
\def\P{\mathbb P}

\def\calb{{\cal B}}

\def\calf{{\cal F}}

\def\calp{{\cal P}}

\def\cals{{\cal S}}

\def\wp{\widehat{P}}
\def\wq{\widehat{Q}}
\def\wu{\widehat{U}}
\def\bx{\bar{X}}
\def\bu{\bar{u}}
\title{Linear-quadratic optimal control under non-Markovian switching}
\author{Fulvia Confortola\thanks{Politecnico di Milano, Dipartimento di Matematica, via Bonardi 9, 20133 Milano, Italy; e-mail: \texttt{fulvia.confortola@polimi.it}}
\and
Marco Fuhrman\thanks{Politecnico di Milano, Dipartimento di Matematica, via Bonardi 9, 20133 Milano, Italy; e-mail: \texttt{marco.fuhrman@polimi.it}}
\and
Giuseppina Guatteri\thanks{Politecnico di Milano, Dipartimento di Matematica, via Bonardi 9, 20133 Milano, Italy; e-mail: \texttt{giuseppina.guatteri@polimi.it}}
\and
Gianmario Tessitore\thanks{Dipartimento di Matematica e Applicazioni, Universit\`a di Milano Bicocca,
via Cozzi 55, 20125 Milano, Italy; e-mail: \texttt{gianmario.tessitore@unimib.it}}
}
\begin{document}
\maketitle
\date{}
\begin{abstract}
We study a finite-dimensional continuous-time optimal control problem
on finite horizon for a controlled diffusion driven by Brownian
motion, in the linear-quadratic case. We admit stochastic coefficients,
possibly depending on an underlying independent marked point process,
so that our model is general
enough to include controlled switching systems
where the switching mechanism is not required to be Markovian.
The problem is solved by means of a Riccati equation,
which a backward stochastic differential equation driven
by the Bronwian motion and by the random measure
associated to the marked point process.
\end{abstract}
\vspace{5mm}
\vspace{5mm}
\section{Introduction}\label{s-Intro}
In order to present and motivate our results
let us consider for a moment  a classical linear-quadratic stochastic
optimal control problem, with a controlled state equation
driven by a $d$-dimensional Brownian motion
$W=(W^1,\ldots,W^d)$ of the form
$$\left\{
    \begin{array}{lll}
    dX_t &=& (A(t)X_t  +B(t) u_t)\, dt + \sum_{j=1}^dC^j(t ) X_t dW^j_t, \\
     X_0 &=& x\in\R^n
    \end{array}
  \right.
$$
and a cost functional
$$
J(u)= \E
\left[\int_0^T (<{S(t)} X_t,X_t> +|u_t|^2)\, dt+ <GX_T,X_T> \right],
$$
where $T>0$ is a fixed finite time horizon,
 $A,B,C^j,S$ are matrix-valued bounded functions
and $S(t)$ and $G$ are non-negative definite. The problem
of minimizing $J(u)$ over all adapted, square-integrable, $\R^k$-valued
processes can be solved via the classical Riccati equation
which provides an optimal feedback control.
More realistic models for many applications require the coefficients
$A,B,C^j,S,G$ to be stochastic. A simple instance is given
by optimization problems for so called {regime-switching}
diffusions, see \cite{Bl-Sw}, \cite{Ji-Chi}, \cite{Ma-Be}, \cite{Mao} among others,
where the controlled process $X$ is assumed to evolve under
a number of regimes, represented  as the elements of a finite set
$K= \{1, \ldots, m\}$, across
which its behavior  can be markedly different. The system is then
described by another stochastic process $(I_t)_{t\ge 0}$,
with values in $K$, which represents the running regime
and which is often assumed to be piecewise constant,
with random positions $\xi_n$ on random time intervals $[T_n,T_{n+1})$,
where $T_n$ are an increasing sequence of switching times.
The dynamic system of interest is now
$$\left\{
    \begin{array}{lll}
    dX_t &=& (A(t,I_t)X_t  +B(t, I_t) u_t)\, dt +\sum_{j=1}^d C^j(t,I_t ) X_t dW^j_t, \\
     X_0 &=& x,
    \end{array}
  \right.
$$
where
$A,B$ and $C^j$ are bounded functions defined on $[0,T] \times K$,
and a similar modification is performed on the cost functional as well.
For example, in Mathematical finance,  to model the  price of a stock in a
financial market, we can use an equation of the form
$$dS_t=\mu(t,I_t)S_t\, dt+\sigma(t,I_t)S_t \,dW_t,$$
where $S$ represent the stock price, $\mu$ and $\sigma$ the appreciation
and volatility rates, which are modulated by the regime process $I$,
which can be understood as  representing the random
environment, the market trends, an economic regime, a credit (reputation) state
as well as other economic factors.
These models are also called
  controlled hybrid diffusion systems
 or jump linear systems and are the object
of intense study, since they are
fairly general and appropriate for a wide variety
of applications.
For some recent applications in risk theory, financial engineering,
and insurance modeling, we refer the reader to
\cite{DiM-Kab-Run}, \cite{Mo}, \cite{Ro-Sch-Sch-Teu}, \cite{Ya-Yi}, \cite{Zh}
and the references therein. Moreover these models have also been used in
manufacturing, communication theory, signal processing, and
wireless networks; see the many references cited in \cite{Ku}.
In the literature, a standard assumption is
that the process $I$ should be a continuous-time
Markov chain with
state space $K$, characterized by its transition
rates, independent of the
Wiener process $W$. In this case the pair $(X,I)$ is a controlled
Markov process with values in $\R^n\times K$, and
extensions of the standard theory allow to solve the
linear-quadratic optimization problem by means
of a system of Riccati equations, indexed by $i\in K$,
see for instance   Chapter 4 in \cite{CoFrTo},
in particular equation (4.17).
It is the purpose of the present paper to generalize this framework
and consider the case of a general piecewise-constant,
non-Markovian process $I$, independent of $W$. In addition, we will consider more general
  regime sets $K$ which can be possibly infinite (even uncountable).
Thus, in the following, the sequence $(T_n,\xi_n)$ (or equivalently
the process $I$) will only be assumed
to be a marked point process, satisfying a mild technical condition
(Assumption (A) below).
To allow for even greater generality
we will consider a controlled state equation of the form
\bee\left\{
      \begin{array}{ll}
        dX_s = & ( A_s X_s  + B_s u_s)\, ds +\sum_{j=i}^d C_s^j X_s dW^j_s,
        \qquad s\in [t,T]\subset [0,T],\\
        X_t = & x,
      \end{array}
    \right.
\eee
with a quadratic
cost functional
\bee
J(t,x,u)= \E^{  \calf_t}
\left[\int_t^T (<{S_s} X_s,X_s> +|u_s|^2) ds+ <GX_T,X_T> \right]
\eee
where now $A,B,C^j,S$ (respectively, $G$) are matrix-valued bounded stochastic
processes (resp. bounded random variable),
 which are assumed
to be predictable with respect to the filtration
 $(\f_t)_{t \ge  0}$ generated by $I$ and $W$
 (resp. $\calf_T$-measurable).
  $S$ and $G$ are non-negative, as before.
 Correspondingly, the control $u$ will also be
 $(\f_t)$-predictable. The use of a random cost functional is customary
 when dealing with stochastic coefficients, but since
 $\E^{  \calf_0}=\E$ this models generalizes the previous ones when $t=0$.
Our main result states that  the (stochastic) value function has the form
\bee
  \inf_u J(t,x,u)=
<P_t x,x>
\eee
where  $P$ is the unique global solution
 to the following Riccati backward stochastic differential equation:
\bee\label{eq:quadratic-intro}\left\{
      \begin{array}{rll}
        -dP_t&= & \dis (A_t' P_t +P_tA_t + C_t' P_t C_t +C_t' Q_t +Q_t C_t- P_t' B_t B_t' P_t + S_t) \,dt \\
               && \dis - Q_t dW_t  -\int_K U_t \, \tilde{\mu}(dt,dx),\\
        P_T&= & G,
      \end{array}
    \right.
\eee
see Theorems  \ref{th:global exist} and \ref{syn-opt-cpn} for more details.
The unknown in  \eqref{eq:quadratic-intro} is a triple $(P,Q,U)$, where $P$ is a matrix-valued  adapted process with cadlag paths and
 $Q=(Q^1,\ldots,Q^d)  $ is matrix-valued $(\calf_t)$-predictable    processes
 and $U$ is a matrix-valued $(\calf_t)$-predictable
 random field defined on $\Omega \times [0,T]\times K$.
Finally, the optimal control is characterized by the optimal feedback control law
$ {u}_s = -B^{'}_sP_{s-}X_s$.
Hence, we solve completely a linear quadratic stochastic optimal
control problem under non-Markovian switching.
\normalcolor
When only the Brownian motion is present the problem has been widely studied. It was introduced by Bismut in \cite{Bi76} as an open problem and firstly solved by Peng \cite{Peng92}  without control dependent noise (as our case). Only more recently a series of papers \cite{KoTa01}, \cite{KoTa02}, \cite{KoTa03} and eventually \cite{Ta}
have solved the more general case with control dependent case. All these results treat the finite horizon case, in \cite{GuaMas1} and \cite{GuaMas2} there are some extensions to the infinite and ergodic case.
We recall moreover that  a linear quadratic stochastic optimal control with Poisson jumps and Markov switching is discussed  in \cite{Shi-Wu}. Here the authors assume that the constrained stochastic Riccati equation admits a solution and obtain an optimal state feedback control and the
value function.

We believe that our results in the non-Markovian case
can be generalized in several directions, for
instance to the case of control on infinite horizon, both for a discounted
or an ergodic cost functional, and to the more difficult situation
when the control affects the diffusion coefficient
(along the lines of \cite{Ta}, where however the Wiener
process is the only source of randomness)
or even when the controlled equation is driven by some
discontinuous integrator in addition to the Brownian motion.
These extensions are left for future work.
\section{General framework and preliminaries.} This section sets out the notation and
some
assumptions that are supposed to hold in the sequel.
We first describe the noise entering the system.
Let $(\Omega,\calf, \P)$ be a complete probability space, where
a standard $d$-dimensional Brownian motion
$W=(W^1,\ldots,W^d)$ is defined as well as an independent
 multivariate point process
(also called marked point process) on a space $ K$.
Next we recall some basic properties of such point processes for which we refer
to \cite{Ja} or
\cite{Bre} or \cite{BrLa}.
We  suppose that   $K$ a Borel space,
i.e. a topological space homeomorphic to a Borel subset of a compact
metric space (some authors call this a Lusin space); thus,
$K$ can be any complete separable metric space.
 The Borel $\sigma$-algebra of $K$ is denoted by $\calb(K)$
 (a similar notation will also be used for other topological spaces as well).
A marked point process
is a double sequence $(T_n,\xi_n)_{n\ge 1}$ such that
the random variables  $T_n$ take values
in $(0,\infty]$ and satisfy
$ T_n<T_{n+1}$ whenever $T_n<\infty$, and
the random variables $\xi_n$ (called \emph{marks}) take values in
$K$ and satisfy
$\xi_n=\De$ whenever
$T_n=\infty$, where $\De$ is a distinguished point in $K$.
We will impose conditions implying that the process is non-explosive, that is
 $T_n\to\infty$.
Let $(\mathcal{F}_t)_{ t \ge 0}$ denote the smallest complete
right continuous filtration generated by $I$ and $W$. Throughout the
paper we only use this filtration. We denote the conditional
expectation with respect to $\calf_t$ by the symbol
$\E^{\calf_t}(\cdot)$ (rather than $\E[\cdot\mid{\calf_t}]$).
We let
$\p$ denote the predictable $\si$-algebra corresponding
to  $(\mathcal{F}_t)_{ t \ge 0}$. By abuse of notation,
we use the same symbol to denote the trace of $\p$ on any subset
 $\Om\times J$ for any interval $J\subset [0,\infty)$. For any auxiliary
measurable space $(G,\g)$, a function on the product $\Om\times J\times G$
which is measurable with respect to $\p\otimes\g$ is also called
{\em predictable}.
To the
 marked point process we can associate a $K$-valued piecewise constant process
 $I$ defined by
 $ I_t=\xi_n$ for $t\in [T_nT_{n+1})$ (and $I_t=k_0$, some given
 point in $K$, for $t\in [0,T_1)$)
  and a random measure $\mu$  on
$((0,\infty) \times K,\mathcal{B}((0,\infty)\times K))$
given by
\bee\label{S1}
\mu(dt,dx)=\sum_{\,T_n<\infty }\delta_{(T_n,\xi_n)}(dt,dx).
\eee
We need the concept of
  compensator (or dual predictable projection) of $\mu$  under $\P$,
relative to the filtration $(\f_t)$.
This is a predictable random measure $((0,\infty) \times K,\mathcal{B}((0,\infty)\times K))$,
denoted $\nu(dt\,dx)$,
satisfying
\bee \label{phi}
\E \int_{0}^\infty \int_K H_t(x)\; \mu(dt\, ,dx)=\E \int_{0}^\infty \int_K H_t(x)\;\nu(dt\,dx),
\eee
for every nonnegative predictable process $H$.
The
measure $\nu$ admits the disintegration:
\bee\label{S2}
\nu(\om,dt,dx)~=~da_t(\om)\,\phi_{\om,t}(dx),
\eee
where $a$ is an increasing c\`adl\`ag predictable process
starting at $a_0=0$ (which is also the  compensator of the
univariate point process
$\mu((0,t]\times K)$, $t\ge 0$)
and
$\phi$ is a transition probability from $(\Om\times(0,\infty),\p)$ into
$(K,\mathcal{K})$
$ $
\noindent We make the following

$ $

\noindent \textbf{Assumption (A)} $\P$-a.s., the process $(a_t)_{t>0}$ has continuous
trajectories.

$ $

It can be proved that Assumption (A) implies that the process
is non-explosive, and in fact it
is equivalent to the requirement that the
jump times $T_n$ are non exposive and totally inaccessible. (A) holds
if and only if, $\P$-a.s.,  $\nu(\{t\}\times K)=0$ for every $t> 0$.
We finally note that we will be interested in a control
problem formulated for a fixed   deterministic time horizon  $T \in (0, \infty)$,
so that we only need to have $W$ defined on $[0,T]$
and $\mu$ a random measure defined on $(0,T]\times K$.
For  any Euclidean space  $E$,  we denote
 by  $<\cdot,\cdot>$   the scalar product  and by $\mathcal{B}(E)$ the Borel $\sigma$-algebra.
  We denote by $\mathcal{S}_n$ the space of symmetric matrices of dimension $n\times n$, and  by $\mathcal{S}^+_n$ its subset of non-negative definite matrices. We denote
by the same symbol $|\cdot|$ both the norm of a vector and the matrix operator norm.
Let $a,b$ be real numbers, $0\le a<b\le T$.
The following classes of processes will be used in the paper.
 \begin{itemize}
\item $ L^p_{\p}(\Omega \times [a,b]; E)$,  for
$p \in [1, \infty]$ denotes the standard $L^p$ space
constructed on the measurable space $(\Omega \times [a,b],\calp)$
endowed with the product measure $\P(d\omega)\,dt$.
It is endowed
with the natural norm
$$|Y |^p_{L^p_\p(\Omega \times [a, b];E)} =\E \int_a^b |Y_s|^p ds $$
for $p<\infty$, replaced by the essential supremum of $|Y|$ for $p=\infty$.
 Elements of this space are identified up to almost sure
 equality with respect to $\P(d\omega)\,dt$.
\item $L^p_{\p}(\Omega;D([a, b];E))$, for $p\in [1,\infty]$,
 denotes the space of adapted processes $Y$ with c\`adl\`ag
paths in $E$ (i.e., right-continuous on $[a,b)$ having finite
left limits on $(a,b]$) such that the norm
$$
\begin{array}{lll}|Y |^p_{L^p_{\p}(\Omega;D([a, b];E))} = \E \sup_{t\in [a,b ]} |Y_t|^p
& \mbox{if}& p<\infty,
\\
|Y |_{L^\infty_{\p}(\Omega;D([a, b];E))} = \mathop{ess\;sup}_{\omega\in\Omega}
\sup_{t\in [a,b ]} |Y_t(\omega)|&
  \mbox{if}& p=\infty
  \end{array}
$$
is finite.
Elements of this space are identified up to indistinguishability.

\begin{remark}{\em The previous notation is justified from the fact that a process $\tilde Y \in L^p_{\p}(\Omega;D([a, b];E))$ is progressively measurable and it is well known that given such a process, it is possible to find
$Y\in L^p_{\p}(\Omega \times [0,T];  \R^{ k})$ such that
$Y=\tilde Y$ $\P(d\omega)\,dt$-a.s.,
}
\end{remark}

\item $L^p_{\p}(\Omega;C([a, b];E))$, for $p\in [1,\infty)$,
 denotes the subspace of
 $L^p_{\p}(\Omega;D([a, b];E))$ consisting of processes with continuous
 paths. It is endowed with the same norm and its elements
 are
 predictable processes.

\item $L^p(a,b,\nu)$, for $ p\in [1, \infty)$,
denotes the set of equivalence classes, with respect  to the measure $\phi_t(\omega,dx)da_t(\omega)\P(d\omega)$, of mappings
$H : \Omega\times (a, b] \times K \rightarrow \mathcal{S}_n$ which are predictable (i.e. $\mathcal{P}\otimes \mathcal{B}(K)$-measurable)
and such that
$$ |H|_{L^p(a,b,\nu)}^p= \E \int_{(a,b]}  \int_K |H_t(x)|^p\; \mu(dt,dx)= \E \int_{(a,b]}
\int_K |H_t(x)|^p \, \nu(dt,dx)  < \infty.$$
\end{itemize}
Moreover we denote with $L^p(\Omega, \f_T,\P;E)$  the subset of
$\P$-equivalence classes of  $L^p(\Omega, \f,\P;E)$ which have
an  $\f_T$-measurable representative, endowed
with the same norm ($p\in [1,\infty]$).
We recall that for any predictable
real function satisfying
   $
\int_0^T \int_K |H_t(y)| \;\phi_t(dy)\, da_t <\infty,
\P$-a.s.
one can define the stochastic integral with respect to $\tilde{\mu}= \mu- \nu$   as the difference of ordinary
integrals with respect to $ \mu$ and $\nu$.
Given an element $H$ of $L^1(0,T,\nu)$, its
 stochastic integral   with respect to $\tilde{\mu}$
  turns out to be a  finite variation martingale on $[0,T]$.
 Moreover if $H$ is in $L^2(0,T,\nu)$ then
 its  stochastic integral with respect to $\tilde{\mu}$  is a square integrable, purely discontinuous martingale with predictable quadratic variation
 $\int_{0}^\cdot \int_K  |H_t(x)|^2\;\phi_t(dx)\, da_t  $.
Finally we recall that the weak property of predictable representation
holds with respect to $(\f_t)$ and $\P$ (see \cite[Example 2.1 (2)]{Bech}).
This means that every square integrable martingale $M$ has a representation
$$M_t =M_0 +\int_0^t Z_s\, dW_s + \int_0^t \int_K U(s,x) \tilde{\mu}(ds,dx)$$
   where $Z\in L^2_{\p}(\Omega \times [0,T]; \R^d)$
   and $U\in L^2(0,T,\nu)$.
\section{Assumptions and statement of the problem}
Throughout the paper we assume that a
probability space
$(\Omega,\calf, \P)$, a Brownian motion
$W$ and an  independent
 multivariate point process $(T_n,\xi_n)_{n\ge 1}$   on a space $ K$ are given,
 satisfying the assumptions in the previous section, in particular
 Assumption (A) that will be recalled in the   statements of the main results.
We consider the following stochastic differential equation
\bee\left\{
      \begin{array}{rll}\label{eq:stato}
        dX_t &= & ( A_t X_t dt + B_t u_t) \,dt + C_t X_t \,dW_t,\\
        X_s &= & x,
      \end{array}
    \right.
\eee
where the unknown process $X$ is $\R^n$-valued
and represents the state of a controlled system, $u$ is the control process
and the initial condition $x\in\R^n$ is deterministic.
 A precise notion of solution to
the state equation (\ref{eq:stato})  is given below.
To stress its dependence
on $u$, $t$, and $x$ we will denote it by $X^{t,x,u}$
when needed.
We introduce a cost functional
of the form
$$
J(t,x,u)= \E^{\calf_t}
\left[\int_t^T (< {S_s} X_s,X_s> +|u_s|^2)\, ds+ <GX_T,X_T> \right]
$$
and we aim at finding an optimal control, relatively
to the given data $(t,x)$, that is
$\bar{u} \in L^2_{\calp}(\Om \times [t, T];\R^k)$ such that
$$J(t, x, \bar{u}) =
\mathop{\rm ess \; inf}_{u \in L^2_{\calp}(\Om \times [t, T];\R^k)} J(t, x, u).$$
We also look for a characterization of the (random) value  function,
that is the essential infimum above.
Elements of  the space $L^2_{\p}(\Omega \times [t,T];  \R^{ k})$
are called admissible controls.
\begin{remark}{\em
The minimization could also be equivalently performed
over $L^2_{\p}(\Omega \times [0,T];  \R^{ k})$ since the values
of the process $u$ over $[0,t]$ are irrelevant.

Another possible formulation consists in  considering control processes $\tilde u$
satisfying $\E\int_0^T|\tilde u_s|^2ds<\infty$
which are only progressively measurable (rather than predictable).
However, given such a process $\tilde u$, it is possible to find
$u\in L^2_{\p}(\Omega \times [0,T];  \R^{ k})$ such that
$u=\tilde u$ $\P(d\omega)\,dt$-a.s., so that the corresponding
trajectories coincide and we clearly have
$J(t,x,u)=J(t,x,\tilde u)$. Therefore the two optimization problems
are essentially the same. If one prefers to use progressively measurable
control processes the optimal feedback law \eqref{contr-ottimo} simplifies
to $\bar{u}_s = -B^{'}_sP_{s}\bx_s$.
}
\end{remark}
We will work under the following general assumptions
on the coefficients.
\begin{hypothesis}\label{hyp:eq stato}
$ $
\begin{itemize}
\item[($A_1$)] We assume that the processes $A,B,C=(C^1,\ldots,C^d)$ satisfy
$$
A \in  L^{\infty}_{\p}(\Omega \times [0,T];  \R^{n \times n}),
\quad
B \in L^{\infty}_{\p}(\Omega \times [0,T];  \R^{n \times k}),
\quad
C^j \in  L^{\infty}_{\p}(\Omega \times [0,T]; \R^{n \times n}),
$$
for $j =1,\ldots, d$.
\item [($A_2$)] $G \in L^{\infty}(\Om,\calf_T, \P; \cals^+_n)$.
\item [($A_3$)] $S \in  L^{\infty}_\p(\Om \times [0,T]; \cals^+_n)$.
\end{itemize}
We denote $M_A,M_B, M_C, M_G, M_S$ nonnegative constants such that
$|G(\om)|\le M_G$ $\P(d\omega)$-a.s. and
$$
|A_t( \om)| \leq  M_A,\quad
|B_t( \om)| \leq  M_B, \quad
|C^j_t( \om)| \leq  M_C, \quad
|S_t( \om)| \leq  M_S, \quad
$$
$\P(d\omega)dt$-a.s. for $j =1,\ldots, d$.
\end{hypothesis}
Next we present  precise statements that ensure that the formulation of the
optimization problem makes sense.
\begin{definition} Given $x \in \R^n$, $t\in [0,T]$
and $u \in L^2_{\p}(\Omega \times [0,T];  \R^{ k})$, a solution to (\ref{eq:stato})
is a process $X \in L^2_{\p}(\Omega;C( [0,T];  \R^{ n}))$
 such that, $\P$-a.s.,
\bee
X_s= x + \int_t^s ( A_rX_r + B_r u_r)\, dr + \sum_{j=1}^d\int_t^s C_r^j X_r \,dW^j_r,
\qquad  s\in [t, T].
\eee
\end{definition}
The following existence and uniqueness result is standard (see \cite{GS},\cite{IW} or \cite{Pro}).
\begin{theorem}\label{th: stima su X} Let assumption ($A_1$) be satisfied. For any $p \geq 2$,
given any $t\in [0,T]$, $x \in \R^n$ and  predictable control  $u$ with
$$\E\left( \int_t^T |u_s|^2 ds\right)^{p/2} < \infty,$$
the equation (\ref{eq:stato}) has a unique  solution $X \in L^p_{\p}(\Omega; C([t, T]; \R^n))$ and
it satisfied the estimate
\bee \label{stima su X}
 \E^{\mathcal{F}_t} \sup_{sÅ\in [t,T ]} |X_s|^p \leq  C_p \left[|x|^p +
 \E^{\mathcal{F}_t} \left(\int^T_t |u_s|^2 ds\right)^{p/2} \right]
\eee
for a suitable constant $C_p$ depending on $p,\,T,\,M_A,\,M_B $ and $M_C$. Notice that $C_p\geq 1$.
\end{theorem}
\section{Solution of the optimal control problem}
\subsection{The  Lyapunov equation.} We start from the linear part of the Riccati
equation. Namely we consider the Lyapunov equation
\bee\label{eq:Lyapunov}\left\{
      \begin{array}{rll}
        -dP_t&= & \dis (A_t' P_t +P_tA_t + C_t' P_t C_t +C_t' Q_t +Q_t C_t+ L_t) \,dt \\
               & &\dis - Q_t dW_t  -\int_K U_t \, \tilde{\mu}(dt,dx),\\
        P_T&= & H,
      \end{array}
    \right.
\eee
where $L \in L^2_{\p}(\Omega \times [0,T];\mathcal{S}_n)$ and $H\in L^2(\Omega, \mathcal{F}_T, \P; \mathcal{S}_n)$. We   use  the shortened notation
\begin{equation}\label{shortened}
    C_t' P_t C_t +C_t' Q_t +Q_t C_t=
    \sum_{j=1}^d[    (C_t^j)' P_t C^j_t +(C_t^j)' Q_t^j +Q_t^j C^j_t]
    ,\qquad
    Q_t dW_t =\sum_{j=1}^dQ_t^j dW^j_t.
\end{equation}
\begin{definition} \label{def:sol-Lyap} A solution to problem (\ref{eq:Lyapunov}) is a process
$(P,Q,U) \in L^2_{\p}(\Omega;D([0, T];  \mathcal{S}_n)) \times L^2_{\p}(\Omega \times [0,T];(\mathcal{S}_n)^d) \times L^2(0,T,\nu)$
that verifies, $\P$-a.s.,
\bee
P_t= H+ \int_t^T (A_s' P_s +P_sA_s + C_s' P_s C_s +C_s' Q_s +Q_s C_s+ L_s) \,ds
-\int_t^T Q_s dW_s  -\int_t^T\int_K U_t \, \tilde{\mu}(ds,dx),
\; t \in [0, T].
\eee
\end{definition}

Proposition \ref{Prop:Lyap-existence-uniqueness} ensures existence and uniqueness of the solution of the Lyapunov equation \eqref{eq:Lyapunov}. We remark that Assumption (A)
is used at this point, but it is not needed in the sequel.
\begin{proposition}\label{Prop:Lyap-existence-uniqueness} Assume Hypotheses ($A_1$).
Then for any   $H \in L^2(\Omega, \mathcal{F}_T, \P; \mathcal{S}_n)$
and $L \in L^2_{\p}(\Omega \times [0,T];\mathcal{S}_n)$
 problem (\ref{eq:Lyapunov}) has a unique  solution $(P,Q,U)$ and we have moreover
\bee \label{stima sol Lyap} \E \sup_{s \in [t,T]}|P_s|^2 + \E \int_t^T |Q_s|^2 ds + \E \int_t^T \int_K |U_s(x)|^2 \nu(ds,dx) \leq C_0\; \E \left[ |H|^2 + \int^T_t|L_s|^2 ds \right],
\eee
for every $t\in [0,T]$ and for some constant $C_0$ depending only on
$T,M_A,M_C$ and the underlying marked point process.
\end{proposition}
\textbf{Proof.}   The proof of this and other similar results
relies on the weak property of predictable representation
  mentioned above.
  In the case of a Poisson
 random measure (possibly however with infinite activity)
 the result was proved in Lemma 2.4 of \cite{Li-Tang}, in Theorem 2.1 in \cite{Ba-Bu-Pa} and in Theorem 53.1 in \cite{Pardoux}.
  The result is also proved in \cite{Bech} in the setting of a nonhomogeneous compensator $ \nu$ assumed to be absolutely
 continuous with respect to Lebesgue measure.
Under the Assumption (A), it is straightforward to generalize the established fixed
point method of proof to the present setting, see for instance
\cite{Co-Fu}.
For this reason we omit the proof and  leave the details to the reader.
\qed

The following result is a key step towards the fundamental relation
(see Proposition \ref{fund-rel+pos+aprioriest}-1).
\begin{theorem} Assume Hypotheses ($A_1$). Let $H \in L^\infty(\Omega, \mathcal{F}_T, \P; \mathcal{S}_n)$,
$L \in L^{\infty}_\p (\Omega \times [0,T];\mathcal{S}_n)$ and let $(P,Q,U)$ be the
unique  solution to (\ref{eq:Lyapunov}).  Then for all
$t\in [0, T], x \in \R^n, u \in L^2_{\p}(\Omega \times [0,T];  \R^{ k})$,
denoting  by  $X^{t,x,u}$ the corresponding solution to (\ref{eq:stato}),
it holds  that, $\P$-a.s.,
\bee \label{Lyap-fund-rel}
<P_t x, x> = \E^{\calf_t }<H X^{t,x,u}_T, X^{t,x,u}_T>
+ \E^{\calf_t } \int_t^T [ <L_s  X^{t,x,u}_s, X^{t,x,u}_s>
- 2< P_sB_s u_s, X^{t,x,u}_s>] ds
\eee
Moreover, for all $t \in [0,T]$,
\bee \label{|P|-extim}
|P_t| \leq C_2 \left[|H|_{L^{\infty}(\Om,\f_T,\P;\cals_n)}
+ (T-t)|L|_{L^{\infty}_\p (\Omega \times [t,T];\mathcal{S}_n)}  \right], \quad \P\mbox{-a.s.}
\eee
where $C_2\geq 1$ is the constant in \eqref{stima su X}.
In particular, we have $P\in L^{\infty}_\p (\Omega \times [t,T];\mathcal{S}_n)$.
\end{theorem}
\textbf{Proof.} \emph{First step.} We first prove  \eqref{Lyap-fund-rel} for $u \in L^8_{\p}(\Om \times [0,T];\R^k)$. The corresponding process
 $X=X^{t,x,u}$   solution to \eqref{eq:stato} then belongs to  $ L^8_{\p}(\Om \times [0,T];\R^n)$ by Theorem \ref{th: stima su X}.
Differentiating by the It\^{o} rule (see e.g. \cite{He-Wang-Yan}, Theorem 9.35) we obtain
$$
\begin{array}{lll}
\dis d <P_s X_s,X_s>  &=& \sum_{i=1}^d [<Q_s^i X_s,X_s> +2<P_s X_s ,C_s^i X_s>]\,dW_s^i
\\
\dis &  & + \int_K<U_s(y)X_s,X_s>\, \tilde{\mu}(ds,dy)
 -[<L_s X_s,X_s> -2<P_s B_s u_s, X_s>]ds.
\end{array}
$$
In order to prove that the local martingale terms have zero mean we introduce an approximating procedure.
Let $\Psi \in C^2(\mathbb{R}^n)$ with $\Psi(y) = 1$ for $|y| \leq 1$, $\Psi(y) = 0$ for $|y| \geq 2$ and $\Psi(y)\in [0, 1], \, \forall y \in \mathbb{R}^n$. Again by the It\^{o} rule
we obtain,  for all integer $N \ge 1$,
\bee \label{5.9}
\begin{array}{lll}
\dis d[\Psi(X_s/N)<P_s X_s,X_s>] &=& N^{-1}F_N(s)ds + G_N(s)dW_s + \Psi(\frac{X_s}{N})\int_K<U_s(y)X_s,X_s>\, \tilde{\mu}(ds,dy)\\
\dis &  & -\Psi(X_s/N)[<L_s X_s,X_s> -2<P_s B_s u_s, X_s>]ds,
\end{array}
\eee
where
$$
\begin{array}{lll}
\dis F_N(s) & = & <\Psi' (\frac{X_s}{N}), [A_s X_s +B_s u_s]><P_s X_s, X_s> \\
\dis & & +2 \sum_{i=1}^d<\Psi' (\frac{X_s}{N}), C^i_s X_s ><P_s C^i_s X_s, X_s>
\\&&
+\frac{1}{2N} \sum_{i=1}^d<\Psi^{''} (\frac{X_s}{N}) C^i_s X_s , C^i_s X_s><P_s X_s, X_s>\\
\dis & &Å+\sum_{i=1}^d<\Psi' (\frac{X_s}{N}), C^i_s X_s ><Q^i_s X_s, X_s>
\end{array}
$$
and for $i=1,...d$
$$
\dis G_N^i(s) = \frac{1}{N}<\Psi' (\frac{X_s}{N}),C^i_s X_s><P_s X_s, X_s>
+\Psi(\frac{X_s}{N})\left(2<P_s C^i_s X_s,X_s> + <Q^i_s X_s, X_s>\right).
$$
It can be easily verified that $\sup_N\E\int^T_t
|F_N(s)|ds <\infty$.

Moreover,
since $\Psi(N^{-1}y) = 0$ and $\Psi'(N^{-1}y) = 0$ if $|y| > 2N$ we have, for all
fixed $N \ge 1$,
$$  \sum_{i=1}^d \E \int_t^T |G_N^i(s)|^2 ds
\leq c N^4 \left(M^2_C T\, \E  \sup_{s\in[t,T]} |P_s|^2
+ \E\int_t^T |Q_s|^2 ds \right) < \infty,$$
 and
$$\E \int_t^T\int_K |\Psi(\frac{X_s}{N})<U_s(y)X_s,X_s>| \, \nu(ds,dy) \leq c N^2 \left(\E \int_t^T\int_K |U_s(y)|^2 \, \nu(ds,dy)\right)^{1/2} < \infty,$$
for a suitable positive constant $c$.

 Finally $<L X,X>$ and $<PBu,X>$ belong to $L^1_{\mathcal{P}}(\Omega \times [t,T];\R)$,
$<P_T X_T,X_T>$ belongs to $L^1(\Omega,\calf_T,\P;\R)$,
and $\Psi(X_s/N)$ boundedly converges to $1$ $\P$-a.s. for all $s$.

Thus, first integrating in $[t, T]$ and then computing conditional expectation with
respect to $\mathcal{F}_t$, and finally letting $N \rightarrow \infty$,
from  (\ref{5.9}) we deduce:
$$
<P_t x, x> = \E^{\mathcal{F}_t} <P_T X_T, X_T> + \E^{\mathcal{F}_t}\int^T_t [<L_s X_s, X_s>
-2<P_s B_s u_s, X_s>]ds.
$$
\emph{Second step.} We prove estimate \eqref{|P|-extim}.
From the {first step} we know that for all $x \in \R^n $, $\P$-a.s.
\bee
<P_t x, x> = \E^{\mathcal{F}_t} <H X_T^{t,x,0}, X_T^{t,x,0}>
+ \E^{\mathcal{F}_t}\int^T_t <L_s X_s^{t,x,0}, X_s^{t,x,0}> \, ds
\eee
and so
\bee
|<P_t x, x>| \leq |H|_{L^{\infty}(\Om,\f_T,\P;\cals_n)}\E^{\mathcal{F}_t} | X_T^{t,x,0}|^2 + |L|_{L^{\infty}_\p(\Omega \times [0,T];\mathcal{S}_n)}
\int^T_t \E^{\mathcal{F}_t}|X_s^{t,x,0}|^2 ds
\eee
and by estimate \eqref{stima su X} with $u=0$ we have,
 for all $x \in \R^n$ with $|x| \leq 1$,
\bee
|<P_t x, x>| \leq C_2|H|_{L^{\infty}(\Om,\f_T,\P;\cals_n)}
 + C_2 (T-t) |L|_{L^{\infty}_\p (\Omega \times [0,T];\mathcal{S}_n)}, \quad \P\mbox{-a.s.}
\eee
such bound, implies the estimate \eqref{|P|-extim}.

\emph{Third step}. We extend \eqref{Lyap-fund-rel} to all the admissible controls.
For a general $u \in L^2_{\mathcal{P}}(\Omega \times [0, T];\R^k)$ we choose a sequence $u_m$  such that $u_m \rightarrow u$ in $ L^2_{\mathcal{P}}(\Omega \times [0, T];\R^k)$
and each $u_m$ is bounded.
By Theorem \ref{th: stima su X}, $X^{t,x,u_m} Å\rightarrow X^{t,x,u}$ in $L^2_{\mathcal{P}}(\Omega;C([t, T]; \R^n))$ and, by the  {second step}, $P \in  L^{\infty}_{\mathcal{P}}(\Omega \times [0, T]; \cals^n)$.
Equality \eqref{Lyap-fund-rel} holds for $u_m$ and $X^{t,x,u_m}$ and it
is easy to verify that we obtain \eqref{Lyap-fund-rel} for $u$ and $X^{t,x,u}$
letting  $m \rightarrow \infty$. For instance, we may verify that
$$
\begin{array}{lll}
\left| \int^T_t <L_s X^{t,x,u_m}_s, X^{t,x,u_m}_s> - <L_s X^{t,x,u}_s, X^{t,x,u}_s >ds\right| & &\\
\leq \left[\left( \sup_{sÅ\in[t,T ]} |X^{t,x,u_m}_s|^2\right)^{1/2} + \left( \sup_{sÅ\in[t,T]} |X^{t,x,u}_s|^2\right)^{1/2}\right]\cdot& &\\
\cdot
\left( \sup_{sÅ\in[t,T ]} |X^{t,x,u_m}_s- X^{t,x,u}_s|^2\right)^{1/2}
T |L|_{L^{\infty}_\p(\Omega \times [0,T];\mathcal{S}_n)} & &
\end{array}
$$
tends to $0$ in $L^1$. The other terms are treated in a similar way.
\qed
\subsection{Existence and  uniqueness  for the Riccati equation}
 In this section we prove the existence of a
unique  solution for the Riccati equation
\bee\label{eq:Riccati}\left\{
      \begin{array}{ll}
       -dP_t&=  \dis (A_t' P_t +P_tA_t + C_t' P_t C_t +C_t' Q_t +Q_t C_t- P_t' B_t B_t' P_t + S_t) \,dt \\
             & \dis - Q_t dW_t  -\int_K U_t \, \tilde{\mu}(ds,dx)\\
        P_T&=  H
      \end{array}
    \right.
\eee
where $H \in L^{\infty}(\Om,\f_T,\P;\cals_n)$  is  a general final datum  while the other coefficients  are the ones introduced in Assumption (A) and Hypothesis \ref{hyp:eq stato}.
We still use the shortened notation \eqref{shortened}.
The occurrence of a quadratic nonlinear term requires a specific approach
to solve the problem, which is classical
when dealing with the Riccati equation, see for instance \cite{BeDPDeMi} for the classical case and
\cite{Peng92}, Section 5, or \cite{Tes}, when the coefficients are random. First we
will find a local solution and then we will prove some a priori estimate for the solution
to guarantee the existence of a global solution. The method we use to prove the a
priori bound is based on the so-called fundamental relation (see Proposition
\ref{fund-rel+pos+aprioriest} below) and
uses, in an essential way, the control-theoretic interpretation of the Riccati equation.
We give the notion of solution for the equation \eqref{eq:Riccati},
to be compared with   Definition \ref{def:sol-Lyap}.
\begin{definition}\label{def:sol-Riccati}  Fix $T_0 \in [0, T]$. A solution for problem (\ref{eq:Riccati})
on the interval $[T_0, T]$ is a triple $(P,Q,U)$ with
$$P\in L^{\infty}_{\p}(\Om;D([T_0, T];\cals_n)), \quad
Q \in L^2_{\mathcal{P}}(\Omega \times [T_0, T];(\cals_n)^d), \quad
U \in L^2(T_0,T,\nu)$$
such that, $\P$-a.s.,
\bee
\begin{array}{lll}
P_t &=& H +  \int_t^T [A_s' P_s + P_s A_s + C'_sP_sC_s + C'_sQ_s + Q_sC_s+S_s ]\, ds\\
&&-
 \int_t^T Q_s dW_s -\int_t^T \int_K U_s(x) \tilde{\mu}(ds,dx) -\int_t^T P_sB_sB_s'P_s \,ds,
 \qquad  t \in [T_0, T].
\end{array}
\eee
\end{definition}
\begin{proposition}\label{prop:localexistence}
(local existence and uniqueness).
Under Hypotheses \ref{hyp:eq stato},  for every $R   >0$ there exists a
$\delta= \delta (R) \in (0, T]$ such that problem (\ref{eq:Riccati}), with $|H|_{L^\infty} \leq R$, has a  unique solution on the interval $[T - \delta, T]$.
\end{proposition}
\textbf{Proof. } Recall the notation $M_B$, $M_S$ for the constants
introduced in Hypothesis
\ref{hyp:eq stato}. Let $C_p$ and $C_0$ be the constants in
\eqref{stima su X} and \eqref{stima sol Lyap} respectively.
We fix  arbitrarily $r>  C_2 R$ and choose $\delta\in ]0, T]$ satisfying
\begin{equation}\label{delta_contrazione}
C_2[R +  \delta (r^2 M^2_B+M_S)] \leq r, \qquad
4  C_0r^2  M^4_B \delta \leq \frac{1}{2}.
\end{equation}
We define
$$B(r) = \{P \in L^2_\p(\Omega;D([T -\delta, T];\mathcal{S}_n))  : \sup_{t \in[T -\delta,T] }|P_t| \leq r \quad \P\mbox{-a.s.} \}$$
and note that
$B(r)$ is a complete metric space when endowed with the distance of
$L^2_\p(\Omega;D([T -\delta, T];\mathcal{S}_n)) $.
We construct a contraction
map $\Gamma : B(r) \rightarrow B(r)$, letting $\Gamma(P) = \widehat{P}$, where $(\wp,\wq,\wu)$ is
the unique solution to the Lyapunov equation
(\ref{eq:Lyapunov}) on the time interval $[T-\de, T]$
 with $L = S -PBB'P $; that is,
\bee
\begin{array}{lll}
\wp_t = H &+& \int_t^T [A_s' \wp_s + \wp_s A_s +
  C'_s\wp_sC_s + C'_s\wq_s + \wq_sC_s+S_s]\, ds\\
&& -\int_t^T \wq_s \,dW_s -\int_t^T \int_K \widehat{U}_s(x) \tilde{\mu}(ds,dx)
-\int_t^T    P_sB_sB_s'P_s \,ds.
\end{array}
\eee
We
first   check that $\Gamma$
maps $B(r)$ into itself. By Proposition \ref{Prop:Lyap-existence-uniqueness} (applied on
$[T -\de, T]$) we know that $\Gamma(P) \in L^\infty_{\calp}(\Om;D( [T -\de, T];\cals_n))$,
so it is enough to show
that for all $t \in[T -\de, T]$ it holds $|\Gamma(P)_t|\leq r$
$\P$-a.s. Thanks to \eqref{|P|-extim} we have, for all $t$,
\bee
\begin{array}{lll}
|\Gamma(P)_t|  &\leq & C_2\left[|G |_{L^{\infty}(\Om,\calf_T, \P;\cals^+_n)} +
{\de}\,|S-P B B' P |_{L^{\infty}_\p (\Omega \times [T-\de ,T];\mathcal{S}_n)}
\right]\\
& \leq &  C_2[R + \de(r^2  M^2_B+M_S)] \leq r, \qquad \P\mbox{-a.s.}
\end{array}
\eee
by \eqref{delta_contrazione}. To check the contraction property,
we take
 $P^1$ and $P^2$ in $B(r)$ and recall (\ref{stima sol Lyap}) obtaining
$$\begin{array}{lll}
\E \sup_{t\in [T - \de, T]} |\Gamma(P^1)_t- \Gamma(P^2)_t|^2
&\leq & C_0\,\E \int_{T-\de}^T |P^1_sB_sB^{'}_s P^1_s -P^2_sB_sB^{'}_sP^2_s|^2 \, ds \\
&\leq & 2C_0\,\E\int_{T-\de}^T
[|(P^1_s-P^2_s)B_sB^{'}_sP^1_s|^2+|P^2_sB_sB^{'}_s(P^1_s-P^2_s)|^2]\, ds   \\
         &\leq & 4 C_0 r^2  M^4_B \delta \;\E \sup_{t \in [T-\de,T]} |P^1_t-P^2_t|^2
\end{array}
$$
so that  $\Gamma$ is indeed a contraction in $B(r)$ by \eqref{delta_contrazione}.

If $P$ is its unique fixed
point, the solution $(P,Q,U)$ of (\ref{eq:Lyapunov}) with $L = S-PBB'P$ is
a solution to (5.11). Notice that
$P\in L^\infty_\p(\Omega;D([T -\delta, T];\mathcal{S}_n)) $ thus $(Q,U)$  are well defined by Proposition 	 \ref{Prop:Lyap-existence-uniqueness}.

Conversely, given two solution $(P^i,Q^i,U^i)$ in $[T-\delta_0 ,T]$, $i=1,2$ let $R'=|P^1|_{L^{\infty}(\Omega;D([T-\delta_0,T];S_n)}+|P^2|_{L^{\infty}(\Omega;D([T-\delta_0,T];S_n)}$ and fix $r'$ and $\delta'\leq \delta_0$ such that $ 4 C_0 r^2  M^4_B \delta'<1/2$ and $C_2[R'+\delta' (r'M_B)^2 +\delta' M_S]\leq r'$ (therefore $r'\geq R'$ since $C_2\geq 1$).

 Both $P^i$ lie in  the ball of radius $r'$ in ${L^{\infty}(\Omega;D([T-\delta',T];S_n)}$ and are fixed points of the above defined mapping $\Gamma$ which is a contraction on such a ball. Therefore they must coincide. Proceeding iteratively we get that $P^1$ and $P^2$ coincide on the whole $[T-\delta_0,T]$.
This implies that the other components
$Q^i,U^i$ must coincide as well by the uniqueness result in Proposition \ref{Prop:Lyap-existence-uniqueness}.
\qed

We prove the following a priori bound for {any} solution with nonnegative final point.
\begin{proposition} \label{fund-rel+pos+aprioriest}
Assume Hypothesis \ref{hyp:eq stato} and let $(P,Q,U)$
be any  solution to (\ref{eq:Riccati}) in the
sense of Definition \ref{def:sol-Riccati} on an interval $[T_0, T]$. Moreover suppose that $H\geq 0$. Then
the following holds.
\begin{enumerate}
\item (The fundamental relation)
For all $t \in[ T_0,T]$, $x \in \R^n$, $u \in L^2_{\calp}(\Om \times [t,T];\R^k)$ it holds
\bee\label{eq:fun-rel}
<P_tx, x>= J(t, x, u)  -\E^{\calf_t}\int^T_t
|u_s + B'_sP_sX^{t,x,u}_s|^2ds, \qquad \P\mbox{-a.s.}
\eee
\item (Positivity) For every $t \in [T_0, T]$ and $x \in \R^n$
we have $<P_t x, x>   \geq 0$ $\P$-a.s.
In particular, $P\in L^\infty_\p(\Omega;D([T_0, T];\cals_n^+))$.
\item (A priori estimate)
For every $t \in [T_0, T]$ we have $|P_t|\leq C_2(|H|_{L^\infty}+ T M_S)$ $\P$-a.s.,
where $C_2$ is the constant in \eqref{stima su X}.
\end{enumerate}

\end{proposition}
\textbf{Proof.} We note that $(P,Q,U)$ is the solution to the Lyapunov equation
\eqref{eq:Lyapunov} with $L=S -PBB'P $. Hence by \eqref{Lyap-fund-rel}
\bee
\begin{array}{lll}
<P_t x, x> &=&\E^{\calf_t }<G X^{t,x,u}_T, X^{t,x,u}_T> +E^{\calf_t } \int_t^T <{S_s}  X^{t,x,u}_s,X^{t,x,u}_s>ds\\
& & -\E^{\calf_t } \int_t^T <P_sB_sB_s'P_s X^{t,x,u}_s,X^{t,x,u}_s>
- \E^{\calf_t }\int_t^T< P_sB_s u_s, X^{t,x,u}_s> ds.
\end{array}
\eee
The  fundamental relation then follows
 adding and subtracting $\E^{\calf_t } \int_t^T|u_s|^2 \,ds$
to the right-hand side.
To prove positivity,
  consider the following closed loop equation, starting at any time
$t \in[ T_0,T]$ with an arbitrary initial data $x\in \R^n$:
\bee\left\{
      \begin{array}{lll}
      d \bx_s &=&[A\bx_s - B_sB_s^{'}P_s\bx_s] \,ds + C_s\bx_s dW_s \\
        \bx_t &= & x.
      \end{array}
    \right.
\eee
Such equation fulfills the hypotheses of proposition \ref{th: stima su X}.
Then applying  the fundamental relation (\ref{eq:fun-rel}) to the control $\bu = -B^{'}P\bx$ and  to
$\bx^{t,x,\bu} = \bx$ we get
$<P_tx, x> = J(t,x,\bu) \geq 0$, $\P$-a.s., which proves the
claim.
   Equality (\ref{eq:fun-rel}), with $u= 0$, gives for all $x \in \R^n$ and
all $t Å\in [T_0, T]$,
\bee
\begin{array}{lll}
<P_tx, x> &\le &J(t,x,0)
\\
&=&\E^{\calf_t }<G X^{t,x,0}_T, X^{t,x,0}_T>
+ \E^{\calf_t } \int_t^T < {S_s}  X^{t,x,0}_s, X^{t,x,0}_s>\, ds\\
& \leq & M_G\, \E^{\calf_t }|X^{t,x,0}_T|^2 +M_S\,\int_t^T \E^{\calf_t }|X^{t,x,0}_s|^2\, ds
\end{array}
\eee
and from (\ref{stima su X}) it follows that
$<P_tx, x> \leq C_2\left[M_G + T M_S \right]\,|x|^2$,
which proves the required estimate.
\qed

Now using the
 a priori bound  in
Proposition  \ref{fund-rel+pos+aprioriest}
 we are in a position to extend the  local existence and uniqueness shown in \ref{Prop:Lyap-existence-uniqueness} to the whole $[0,T]$.
\begin{theorem}\label{th:global exist} Suppose that Assumption (A)
 and Hypothesis  \ref{hyp:eq stato} hold true.
 Then the Riccati equation \eqref{eq:Riccati} with $H=G$ has a unique solution
$(P,Q,U$) such that $P\in L^{\infty}_{\calp}(\Om;D( [0, T];\cals^+_n))$,
$Q \in L^2_{\calp}(\Om \times [0, T];(\cals_n)^d)$ and $U \in L^2(0,T,\nu) $.
\end{theorem}
\textbf{Proof.}
Let $R=C_2(M_G+TM_S)$-
By Proposition \ref{prop:localexistence} there exists a unique solution $(P,Q,U$) of equation \eqref{eq:Riccati}  in $[T-{\de}(R),T]$ . Moreover
by Proposition \ref{fund-rel+pos+aprioriest} we know that  $|P_{T-\delta(R)}| \leq R$ and $P_{T-\delta(R)}\geq 0$. Then we can again  can apply the local existence in  $ [T -2{\de}(R), T - {\de}(R)]$ with final datum $P_{T-\delta(R)}$.
 We can then argue iteratively
and cover the whole interval $[0, T]$ by a finite number of intervals of length $\delta(R)$ and obtain the required global solution.
Uniqueness is proved in a similar way: we already know that any two solutions
$(P^i,Q^i,U^i)$, $i=1,2$ must satisfy $0\leq P^i_t\leq RI$ for all $t\in [0,T]$. Using iteratively the uniqueness result in Proposition \eqref{prop:localexistence}, fist on $[T-\de(R),T]$ then on $[T-2\de(R),T-\de(R)]$ and so on we get that they must coincide.

Finally, the fact that $P$ takes values in $\cals^+_n$ follows from the positivity property in
Proposition \ref{fund-rel+pos+aprioriest}.
\qed
\subsection{Synthesis of the optimal control}
The following theorem provides a solution to the control problem.
\begin{theorem} \label{syn-opt-cpn} Suppose that Assumption (A)
 and Hypothesis  \ref{hyp:eq stato} hold true. Fix $t\in [0,T]$ and $x \in \R^n$.
 Then the following holds.
\begin{enumerate}
\item  There exists a unique optimal control
$\bar{u} \in L^2_{\calp}(\Om \times [t, T];\R^k)$.
\item If $\bx=X^{t,x,\bu}$ denotes the corresponding solution  (that is the
optimal state), then $\bx$ is the unique solution to the closed loop equation on $[t,T]$:
\bee\label{closed loop}
\left\{
  \begin{array}{rll}
    d \bx_s &= &[A_s\bx  - B_sB_s^{'}P_s\bx_s]\, ds + C_s\,\bx_s dW_s, \\
    \bx_t &= &x.
  \end{array}
\right.
\eee
\item The following feedback law holds $\P$-a.s. for almost every $s \in [t,T]$:
\bee \label{contr-ottimo}
\bar{u}_s = -B^{'}_sP_{s-}\bx_s.
\eee
\item The value function, i.e. the optimal cost,
 is given by $J(t, x, \bar{u}) =  <P_t x, x>$,  $  \P$-a.s.
\end{enumerate}
\end{theorem}
\textbf{Proof.}
The optimal control, if it exists, is unique by the strict convexity of
the map  $u\mapsto J(t, x, u) $ on $ L^2_{\calp}(\Om \times [t, T];\R^k)$.
Let $(P,Q,U)$ be the unique solution to the Riccati equation (\ref{eq:Riccati}).
From the fundamental relation
\eqref{eq:fun-rel}  we have
$$J(t, x, u) =  <P_t x, x> +   \,  \E\int^T_0 |u_s + B^{'}_sP_s \bx^{t,x,u}_s|^2 ds
=  <P_t x, x> +   \,  \E\int^T_0 |u_s + B^{'}_sP_{s-} \bx^{t,x,u}_s|^2 ds,
$$
where the last inequality follows from the fact that $P_t=P_{t-}$, $\P(d\omega)\,dt$-a.s., since $P_t$ has c\`adl\`ag paths.
Then $J(t, x, u) \geq \, <P_tx, x>$ for all $u \in L^2_{\calp}(\Om \times [t, T];\R^k)$
and the equality holds if
and only if (\ref{contr-ottimo}) holds, that is, if and only if $X=\bx$ solves (\ref{closed loop}) and $u =\bar{u}.$
\qed





\end{document}